\documentclass{article}
\usepackage{amsmath,amssymb}
\newtheorem{defn}{Definition}[section]\newtheorem{thm}[defn]{Theorem}

\newtheorem{lem}[defn]{Lemma}\newtheorem{rem}[defn]{Remark}
\begin{document}
\textwidth=125mm\textheight=195mm

\begin{center}{\Large\bf Quantum de Rham complex on ${\mathbb A}_{h}^{1|2}$}\end{center}

\begin{center}{Salih Celik}\footnote{sacelik@yildiz.edu.tr}\end{center}
\begin{center}{Department of Mathematics, Yildiz Technical University, DAVUTPASA-Esenler, Istanbul, 34220 TURKEY.}\end{center}



\begin{abstract}
Introducing $h$- and $h'$-deformations of ${\mathbb Z}_2$-graded (1+2)- and (2+1)-spaces, denoted by ${\mathbb A}_h^{1|2}$ and ${\mathbb A}_{h'}^{2|1}$, a two-parameter first order differential calculus, de Rham complex, on ${\mathbb A}_{h}^{1|2}$ explicitly constructed.
\end{abstract}

\section{Introduction}\label{sec1}
Quantum superspaces as non-commutative spaces are introduced in \cite{manin1} to represent quantum supergroups. They are classified as $q$-deformed and $h$-deformed or non-standard. A non-commutative differential calculus on $q$-deformed spaces was first suggested by Woronowicz \cite{woro}. Using this approach, a bicovariant differential calculus on quantum 3d-superspace was given in \cite{celik3}. The extension of two parameters of this calculus was introduced in \cite{sucelik}.

Another approach to non-commutative differential calculus was developed by Wess and Zumino \cite{wz}. They have proposed the quantum space and its dual space as parts of a differential graded algebra, the quantum de Rham complex of the quantum space. The natural extension to superspace of this case have been studied by many authors(for example, see \cite{celik1,kob-u,soni} and references therein).

De Rham complex contains a first order differential calculus with extension to higher order. A two parameter first order differential calculus on the quantum superplane ${\mathbb A}_h^{1|1}$ that is covariant under the co-action of the quantum supergroup GL$_{h,h'}(1|1)$ is given in \cite{celik2}. In the present paper, the first order differential calculus on the quantum superspace ${\mathbb A}_h^{1|2}$, which is covariant under the co-action of the quantum supergroup GL$_{h,h'}(1|2)$ given in \cite{celik4}, will be introduced by following Wess-Zumino type differential calculus. This calculus is set up with the commutation relations between the generators of function rings $O({\mathbb A}_h^{1|2})$ and their differentials as the generators of $O({\mathbb A}_{h'}^{2|1})$.

In Section 2 we give the definitions of the $(h,h')$-deformed quantum supergroup and of the coordinate rings of $(h,h')$-deformed superspaces. In Section 3 we set up a first order differential calculus on the $h$-deformed superspace ${\mathbb A}_h^{1|2}$. Section 4 concludes with a brief discussion of the contraction procedure for $\Omega({\mathbb A}_h^{1|2})$.

\noindent{\bf Notation.} The grade (or degree) of a homogenous element $a$ belonging to a ${\mathbb Z}_2$-graded vector space (or superspace) is denoted by $\tau(a)=i$, $i\in {\mathbb Z}_2$. If $\tau(a)=0$, then we will call the element $a$ is {\it even} and if $\tau(a)=1$, it is called {\it odd} element.
Throughout the paper, as in classical case, one assumes that odd (Grassmann) objects anticommute among themselves. All our objects are defined over a ground field ${\mathbb K}$ as ${\mathbb R}$ or ${\mathbb C}$ and ${\mathbb K}'$ as the set of Grassmann numbers.

\section{Review of basic structures}\label{sec2}

In this section, we will briefly make mention of the ${\mathbb Z}_2$-graded Hopf algebra ${ O}({\rm GL}_{h,h'}(1|2))$ and $(h,h')$-deformed quantum superspaces from \cite{celik4}, where $h$ and $h'$ are both odd (or Grassmann) numbers.


\subsection{The quantum super Hopf algebra ${ O}({\rm GL}_{h,h'}(1|2))$}\label{sec2.1}

Let $a$, $b$, $c$, $d$, $e$ $\alpha$, $\beta$, $\gamma$ and $\delta$ be elements of a free superalgebra ${ A}$ and we write them in a matrix form as
$T = \begin{pmatrix} a & \alpha & \beta \\ \gamma & b & c \\ \delta & d & e \end{pmatrix} =(t_{ij})$
where latin letters are even and greek letters are odd.

It can be obtained $(h,h')$-deformed commutation relations between the matrix entries of $T$ from the standard FRT construction \cite{FRT}.
There are two things that are important when we find the commutation relations between the matrix elements of $T$. The first is to give the commutativity relations we expect for the coordinates on the classical supergroup ${\rm GL}(1|2)$ formed by 3x3 super matrices when $h=0=h'$. The second is that we need that relations when setting up the differential calculus.

The relations in the following theorem are presented in \cite{celik4}. We consider the $R$-matrix \cite{celik4}
\begin{equation} \label{2.1}
R_{h,h'} = \begin{pmatrix}
1+hh' & 0 & h' & 0 & 0 & 0 & -h' & 0 & 0 \\
0 & 1 & 0 & 0 & 0 & 0 & 0 & 0 & 0 \\
-h & 0 & 1 & 0 & 0 & 0 & hh' & 0 & -h' \\
0 & 0 & 0 & 1 & 0 & 0 & 0 & 0 & 0 \\
0 & 0 & 0 & 0 & 1 & 0 & 0 & 0 & 0 \\
0 & 0 & 0 & 0 & 0 & 1 & 0 & 0 & 0 \\
h & 0 & hh' & 0 & 0 & 0 & 1 & 0 & -h'\\
0 & 0 & 0 & 0 & 0 & 0 & 0 & 1 & 0 \\
0 & 0 & h & 0 & 0 & 0 & h & 0 & 1-hh' \\
\end{pmatrix}
\end{equation}
which satisfies the ${\mathbb Z}_2$-graded Yang-Baxter equation.

In the following theorem we assume that the odd matrix elements of $T$ are anti-commutative with $h$ and $h'$.

\begin{thm} \cite{celik4} 
The matrix entries of $T$ satisfy the $(h,h')$-deformed commutation relations
\begin{align} \label{2.2}
a\alpha &= (1+hh')\alpha a - h'(da + \alpha\delta), \quad a\beta = \beta a + h'(a^2 - ea - \beta\delta) - h\beta^2, \nonumber\\
a\gamma &= (1+hh')\gamma a + h(\gamma\beta - ca), \quad ac = ca - hc\beta - h'\gamma a + hh'\gamma\beta,\nonumber\\
a\delta &= \delta a + h(a^2 -ea + \delta\beta) + h'\delta^2, \quad ad = da + h\alpha a + h'd\delta - hh'\alpha\delta, \nonumber\\
ae &= ea + h\beta(a-e) + h'(e-a)\delta, \quad \alpha\beta=-(1+hh')\beta\alpha + h' (\beta d + e\alpha), \nonumber\\
\alpha\gamma &= -\gamma\alpha, \quad \alpha c = c\alpha, \quad \alpha\delta = -\delta\alpha - h a\alpha + h' \delta d - hh' ad, \nonumber\\
\alpha d &= d\alpha + h' d^2, \quad \alpha e = e\alpha + h \beta\alpha + h'ed - hh' d\beta, \nonumber\\
\beta\gamma &= - \gamma\beta + h c\beta - h'\gamma a - hh'ca, \quad \beta c = (1-hh')c\beta - h'(ca + \gamma\beta), \nonumber\\
\beta\delta &= - \delta\beta + (h\beta + h'\delta)(e - a), \quad \beta d = d\beta + h\alpha\beta + h'de - hh' e\alpha, \nonumber\\
\beta e &= e\beta + h'(e^2 - ea - \delta\beta) - h \beta^2, \quad \gamma c = c\gamma + hc^2, \\
\gamma\delta &= - (1+hh')\delta\gamma + h (e\gamma + \delta c), \quad \gamma d = d\gamma, \nonumber\\
\gamma e &= e\gamma + h ec - h' \delta\gamma - hh' c\delta, \quad c\delta = \delta c - h ec - h'\delta\gamma - hh' \gamma e, \nonumber\\
cd &= dc, \quad ce = (1-hh')ec + h'(e\gamma - \delta c), \,\,\, \delta d = (1-hh')d\delta + h(\alpha\delta - da), \nonumber\\
\delta e &= e\delta + h(e^2 - ea + \beta\delta) + h' \delta^2, \quad de = (1-hh')ed + h(\beta d - e\alpha), \nonumber\\
\alpha^2 &= h' \alpha d, \quad \beta^2 = h'\beta (e-a), \quad \gamma^2 = h\gamma c, \quad \delta^2 = h\delta (e-a), \nonumber\\
b \,t_{ij} &= t_{ij} b, \quad a(h\beta + h'\delta) = (h\beta + h'\delta)a, \quad e(h\beta + h'\delta) = (h\beta + h'\delta)e, \nonumber
\end{align}
if and only if
\begin{equation*}
R_{h,h'} T_1T_2 = T_2T_1 R_{h,h'}
\end{equation*}
where $T_1=T\otimes I_3$ and $T_2=PT_1P$ and $P$ is the super permutation matrix.
\end{thm}
\noindent The classical limit is obtained by setting $h$ and $h'$ equal to zero.

\begin{defn} 
The superalgebra ${ O}({\rm M}_{h,h'}(1|2))$ is the quotient of the free algebra ${\mathbb K}\{a,b,c,d,e,\alpha,\beta,\gamma,\delta\}$ by the two-sided ideal $J_{h,h'}$ generated by the relations $(\ref{2.2})$ of Theorem 2.1.
\end{defn}

A quantum supergroup (super-Hopf algebra) can be regarded as a generalization of the notion of a supergroup. It is defined by
\begin{equation*}
{ O}({\rm GL}_{h,h'}(1|2)) = { O}({\rm M}_{h,h'}(1|2))[t]/(t \,{\rm sdet}-1).
\end{equation*}
Here the quantum superdeterminant for the supermatrix $T$ is given by
\begin{equation*}
{\rm sdet}(T) = {\rm sdet}(A) \cdot (D-CA^{-1}B)^{-1}
\end{equation*}
where
\begin{equation*}
A = \begin{pmatrix} a & \alpha \\ \gamma & b \end{pmatrix}, \quad B = \begin{pmatrix} \beta \\ c \end{pmatrix}, \quad
C = \begin{pmatrix} \delta \\ d \end{pmatrix}^t, \quad D = \begin{pmatrix} e\end{pmatrix}.
\end{equation*}

\subsection{The $(h,h')$-deformed quantum superspaces}\label{sec2.2}

In this subsection, we introduce an $h$- and $h'$-deformations of the superspaces ${\mathbb A}^{1|2}$ and ${\mathbb A}^{2|1}$ given in \cite{celik4}.

Let ${\mathbb K}\langle x,\theta_1,\theta_2 \rangle$ be a free algebra with unit generated by $x$, $\theta_1$ and $\theta_2$, where the coordinate $x$ is even, the coordinates $\theta_1$ and $\theta_2$ are odd.

\begin{defn} \cite{celik4} 
Let $I_h$ be the two-sided ideal of ${\mathbb K}\langle x,\theta_1,\theta_2 \rangle$ generated by the elements $x\theta_1-\theta_1 x$,
$x\theta_2-\theta_2 x-hx^2$, $\theta_1\theta_2+\theta_2\theta_1$, $\theta_1^2$ and $\theta_2^2+h\theta_2 x$. The quantum superspace ${\mathbb A}^{1|2}_h$ with the function algebra
$${ O}({\mathbb A}^{1|2}_h) = {\mathbb K}\langle x,\theta_1,\theta_2 \rangle/I_h$$
where the parameter $h$ $(h\ne0)$ is a Grassmann number $(h^2=0)$, is called ${\mathbb Z}_2$-graded quantum space (or quantum superspace).
We call ${ O}({\mathbb A}_h^{1|2})$ the algebra of functions on the
${\mathbb Z}_2$-graded quantum space ${\mathbb A}_h^{1|2}$.
\end{defn}

\noindent In accordance with this definition, we have
\begin{equation} \label{2.3}
x\theta_1 = \theta_1 x, \,\,\, x \theta_2 = \theta_2 x + h x^2, \,\,\, \theta_1 \theta_2=-\theta_2\theta_1, \,\,\, \theta_1^2=0, \,\,\ \theta_2^2=-h\theta_2 x.
\end{equation}

\noindent{\bf Note 1.} Although, in the $h$-deformed case, the square of $\theta_2$ is not zero, we have
\begin{equation*}
(\theta_2)^3 = h \,\theta_2^2 \,x = -h^2 \,\theta_2 \,x^2 = 0.
\end{equation*}

\begin{defn} 
Let ${ O}({\mathbb A}_{h'}^{2|1}):=\Lambda({\mathbb A}_h^{1|2})$ be the algebra with the odd generator $\varphi$ and even generators $y_1$, $y_2$ satisfying the relations
\begin{equation} \label{2.4}
\varphi^2 = -h' y_2\varphi, \quad \varphi y_1 = y_1 \varphi, \quad \varphi y_2 = y_2 \varphi - h' y_2^2, \quad y_1y_2 = y_2y_1
\end{equation}
where the parameter $h'$ ($h'\ne0$) is a Grassmann number. We call ${ O}({\mathbb A}_{h'}^{2|1})$ the exterior algebra of the quantum superspace ${\mathbb A}_h^{1|2}$.
\end{defn}

\begin{rem} 
Although the commutation relations between the generators of exterior algebra given in the above definition seem to be obtained by taking $-h'$ instead of $h'$ in \cite{celik4}, due to the odd character of $h'$, the case in reality is not exactly so. An explanation about this case will be given in Sec. 4 (see, Note 5).
\end{rem}

\noindent As it is known, $O({\mathbb A}^{1|2})$ is comodule algebra over the bialgebra $O({\rm M}(1|2))$. The following theorem gives a quantum version of
this fact \cite{celik4}.

\begin{thm} 
The algebras $O({\mathbb A}_{h}^{1|2})$ and $O({\mathbb A}_{h'}^{2|1})$ are both left comodule algebras of the bialgebra $O({\rm M}_{h,h'}(1|2))$ with left coaction $\delta_L$ determined by
\begin{align} \label{2.5}
\delta_L({\bf x}) &= T\dot{\otimes}{\bf x}, \quad {\rm and} \quad \delta_L(\hat{\bf x}) = T\dot{\otimes}\hat{\bf x},
\end{align}
where ${\bf x} = (x, \theta_1, \theta_2)^t$ and $\hat{\bf x} = (\varphi, y_1, y_2)^t$.
\end{thm}

\noindent{\bf Note 2.} Using the $R$-matrix $\hat{R}_{h,h'}$ given in (\ref{2.1}), we can write the relations in (\ref{2.3}) as follows
\begin{equation} \label{2.7}
{\bf x}\otimes{\bf x} = \hat{R}_{h,h'} \,{\bf x}\otimes{\bf x}.
\end{equation}

\begin{rem} 
In $q$-deformation, there is no point in writing the matrix elements of the matrix $\hat{R}_{p,q}$ to the right or left of the tensor product. However, this is very important in $h$-deformation, and care must be taken if it is written on the right hand side.
\end{rem}

\section{The quantum de Rham complex on ${\mathbb A}_h^{1|2}$}\label{sec3}

In this section, we set up de Rham complex, a first order differential calculus, on the superspace ${\mathbb A}_h^{1|2}$. It contains functions on
${\mathbb A}_h^{1|2}$ and their differentials as differential one-forms.

We first note that the quantum supergroup ${\rm GL}_{h,h'}(1|2)$ can be considered as a quantum automorphism supergroup of a pair of noncommutative linear superspaces ${\mathbb A}_h^{1|2}$ and ${\mathbb A}_{h'}^{2|1}$ defined by the relations (\ref{2.3}) and (\ref{2.4}), respectively.

\begin{defn} 
Let ${ A}$ be an arbitrary algebra with unity and $\Omega$ be a bimodule over ${ A}$. A first order ${\mathbb Z}_2$-graded differential calculus over ${ A}$ is a pair $(\Omega,{\sf d})$ where ${\sf d}: { A}\longrightarrow\Omega$ is a linear mapping such that, so called ${\mathbb Z}_2$-graded Leibniz rule,
$${\sf d}(fg) = ({\sf d}f)\, g + (-1)^{\tau(f)} f\, ({\sf d}g) \quad \mbox{for any} \quad f,g\in { A}$$
and $\Omega$ is the linear span of elements of the form $a\cdot{\sf d}b\cdot c$ with $a,b,c\in { A}$.
\end{defn}
\noindent A ${\mathbb Z}_2$-graded differential algebra over $A$ is a ${\mathbb Z}_2$-graded algebra $\Omega=\bigoplus_{n\ge0} \Omega^n$, with the linear map ${\sf d}$ of degree 1 such that ${\sf d}^2:=0$ and ${\mathbb Z}_2$-graded Leibniz rule holds. Here we assume that $\Omega^0:= A$ and $\Omega^{<0} = 0$.

Differential structures are not unique even classically, and even more non-unique in the $q$-deformed case. There is, however, the following theorem (Theorem 3.2) says that differential structure in the $h$-deformed case is unique.

To set up a quantum de Rham complex, the differential calculus, on the quantum superspace ${\mathbb A}_h^{1|2}$ we actually choose the cotangent space or differential 1-forms . Since one can multiply forms by functions from the left and right, this must be an $\Omega$-bimodule.

\subsection{The relations between coordinates and differentials} \label{sec3.1}

Based on Wess and Zumino's suggestion, we can think of the superalgebras $O({\mathbb A}_h^{1|2})$ and $O({\mathbb A}_{h'}^{2|1})$ as parts of a differential graded algebra $\Omega({\mathbb A}_h^{1|2})$, the quantum de Rham complex of the quantum superspace ${\mathbb A}_h^{1|2}$. Thus, we introduce the first order differentials of the generators of ${ O}({\mathbb A}_h^{1|2})$ as ${\sf d}x=\varphi$, ${\sf d}\theta_1=y_1$ and ${\sf d}\theta_2=y_2$. Then the differential ${\sf d}$ is uniquely defined by the conditions in Definition 3.1 and the commutation relations between the differentials have the form
\begin{align} \label{3.1}
{\sf d}x\wedge{\sf d}x &= -h' {\sf d}\theta_2\wedge {\sf d}x, \quad {\sf d}x\wedge {\sf d}\theta_1 = {\sf d}\theta_1\wedge{\sf d}x, \nonumber \\
{\sf d}x\wedge{\sf d}\theta_2 & =  {\sf d}\theta_2\wedge{\sf d}x - h' {\sf d}\theta_2\wedge {\sf d}\theta_2, \quad
   {\sf d}\theta_1\wedge{\sf d}\theta_2 = {\sf d}\theta_2\wedge {\sf d}\theta_1.
\end{align}

In this case, the quantum de Rham complex $\Omega({\mathbb A}_h^{1|2})$ is generated by the elements of the set $\{x,\theta_1,\theta_2,\varphi,y_1,y_2\}$ by adding the nine cross-commutation relations satisfied between the elements of $O({\mathbb A}_h^{1|2})$ and $O({\mathbb A}_{h'}^{2|1})$, which will be given in the following theorem, to the relations (\ref{2.3}) and (\ref{3.1}).

\begin{thm} 
There exists unique left covariant ${\mathbb Z}_2$-graded first order differential calculus $\Omega({\mathbb A}_h^{1|2})$ over the algebra $O({\mathbb A}_h^{1|2})$ with respect to the super-Hopf algebra ${ O}({\rm GL}_{h,h'}(1|2))$ such that $\{{\sf d}x,{\sf d}\theta_1,{\sf d}\theta_2\}$ is a free right
$O({\mathbb A}_h^{1|2})$-module basis of $\Omega({\mathbb A}_h^{1|2})$. The bimodule structure for this calculus is determined by the relations
\begin{align} \label{3.2}
x \,{\sf d}x &= (1 + hh') \,{\sf d}x\, x + h' \,({\sf d}\theta_2 \,x - {\sf d}x \,\theta_2),  \quad  x \,{\sf d}{\theta_1}={\sf d}{\theta_1}\, x, \nonumber \\
x \, {\sf d}\theta_2 &= {\sf d}\theta_2 \, x - h\, {\sf d}x \, x + h' \,{\sf d}\theta_2 \theta_2 + hh' \,{\sf d}x \,\theta_2, \quad
  {\theta_1}\, {\sf d}x = - {\sf d}x \,{\theta_1}, \nonumber\\
{\theta_1} \,{\sf d}{\theta_1} &= {\sf d}{\theta_1} \,{\theta_1},  \quad {\theta_1} \,{\sf d}{\theta_2} = {\sf d}{\theta_2} \,{\theta_1}, \\
{\theta_2} \,{\sf d}x &= - {\sf d}x \,{\theta_2} - h\, {\sf d}x \, x - h' \,{\sf d}\theta_2 \,\theta_2 - hh' \,{\sf d}\theta_2 \,x, \quad
  {\theta_2} \,{\sf d}{\theta_1} = {\sf d}{\theta_1} \,{\theta_2}, \nonumber\\
{\theta_2} \,{\sf d}{\theta_2} &= (1 - hh') \,{\sf d}{\theta_2} \, {\theta_2} - h \,({\sf d}x \,\theta_2 + {\sf d}\theta_2\, x). \nonumber
\end{align}
\end{thm}

\noindent{\it Sketch of proof}. To find nine cross-commutation relations satisfied between the generators of $O({\mathbb A}_h^{1|2})$ and
$O({\mathbb A}_{h'}^{2|1})$, we first write $x{\sf d}\theta_1$, $x{\sf d}\theta_2$, etc. in terms of ${\sf d}\theta_1 x$, ${\sf d}\theta_2 x$, etc. as follows
\begin{align*}
x_i \,{\sf d}x_j &= \sum_{k,l} B_{kl}^{ij} \,{\sf d}x_k \,x_l.
\end{align*}
In contrast to $q$-deformation (which exists forty-one constants there, cf. \cite{celik3}), here we have eighty-one indeterminate constants and we can determine them in four steps as Manin says in \cite{manin2}.

(1) We apply the differential {\sf d} from the left to the relations (\ref{2.3}) and use relations associated with them from cross-commutation relations. In this case, a part of the constants are eliminated.

(2) We apply the operator {\sf d} from the left to cross-commutation relations and compare them with the relations (\ref{3.1}). In this case, more than half of those coefficients are eliminates.

(3) We use compatibility\footnote{Compatibility with the left coaction of $O({\rm GL}_{h,h'}(1|2))$ means that (\ref{2.5}) defines a graded differential algebra homomorphism
$$\Omega({\mathbb A}_h^{1|2})\longrightarrow O({\rm GL}_{h,h'}(1|2))\otimes \Omega({\mathbb A}_h^{1|2})$$}
with the left coaction of $O({\rm GL}_{h,h'}(1|2))$. This leaves one free parameter, say $B_{12}^{21}$, which is not dependent on $h$ and/or $h'$.

(4) The parameter $B_{12}^{21}$ is fixed by checking associativity of the cubics. We find that the parameter $B_{12}^{21}$ should be equal to $1$.
\hfill $\square$

\noindent{\bf Note 3.} The relations in (\ref{3.2}) and (\ref{3.1}) with the matrix $\hat{R}$ can be written in the form
\begin{align} \label{3.3}
(-1)^{\tau(x_i)} \,x_i \,{\sf d}x_j &= \sum_{k,l} (-1)^{\tau(\hat{R}_{kl}^{ij})} \,\hat{R}_{kl}^{ij} \,{\sf d}x_k \,x_l, \\
(-1)^{\tau(x_i)} \,{\sf d}x_i\wedge {\sf d}x_j &= \sum_{k,l} (-1)^{\tau({\sf d}x_k)} \,\hat{R}_{kl}^{ij} \,{\sf d}x_k\wedge {\sf d}x_l,
\end{align}
respectively.

It is possible to define two types of the star operation on Grassmann generators. If $\xi$ and $\eta$ are two Grassmann generators and $c$ is a complex number and $\bar{c}$ its complex conjugate, the star operation, denoted by $\star$, is defined by
\begin{equation*}
(c\xi)^\star = \bar{c}\xi^\star, \quad (\xi\eta)^\star = \eta^\star \xi^\star, \quad (\xi^\star)^\star = \xi
\end{equation*}
and the superstar operation, denoted by $\#$, is defined by
\begin{equation*}
(c\xi)^\# = \bar{c}\xi^\#, \quad (\xi\eta)^\# = \xi^\# \eta^\#, \quad (\xi^\#)^\# = -\xi.
\end{equation*}

We know, from \cite{celik4}, that the algebra $O({\mathbb A}_h^{1|2})$ admits an involution with respect to the deformation parameters, one of them is pure imaginary and the other is real.

\begin{lem} \cite{celik4} 
If $\bar{h}=-h$, then the algebra $O({\mathbb A}_{h}^{1|2})$ supplied with the involution determined by
\begin{equation} \label{3.5}
x^\star = x, \quad \theta_1^\star = \theta_1, \quad \theta_2^\star = \theta_2 - h\,x
\end{equation}
becomes $\star$-algebra.
\end{lem}

The following theorem says that the triple $(\Omega,{\sf d};\star)$ is a $\star$-calculus.

\begin{thm} 
If $\bar{h}=-h$ and $\bar{h}'=h'$, then the first order differential calculus $\Omega({\mathbb A}_h^{1|2})$ is ${\mathbb Z}_2$-graded $\star$-calculus over $O({\mathbb A}_h^{1|2})$.
\end{thm}

\noindent{\it Proof}
We can define the star operation on the differentials as follows \cite{celik4}
\begin{align} \label{3.6}
({\sf d}x)^\star &= {\sf d}x + h' \,{\sf d}\theta_2, \quad ({\sf d}\theta_1)^\star = -{\sf d}\theta_1, \quad ({\sf d}\theta_2)^\star = -{\sf d}\theta_2.
\end{align}
So, it can be shown that the relations (\ref{3.2}) are indeed invariant under the star operation.
\hfill $\square$

\subsection{The relations with partial derivatives}\label{sec3.2}

We will complete the calculus by giving the following two theorems. We first introduce commutation relations between the coordinates of the quantum superspace and their partial derivatives.

Let us begin by introducing the partial derivatives of the generators of the algebra. Since $(\Omega,{\sf d})$ is a left covariant differential calculus, for any element $u$ in ${ O}({\mathbb A}_h^{1|2})$ there are uniquely determined elements $\partial_i(u)\in { O}({\mathbb A}_h^{1|2})$ such that
\begin{equation}\label{3.7}
{\sf d}u = {\sf d}x \partial_x(u) + {\sf d}\theta_1 \partial_{\theta_1}(u) + {\sf d}\theta_2 \partial_{\theta_2}(u).
\end{equation}
For consistency, the degrees of the derivatives $\partial_x$, $\partial_{\theta_1}$ and $\partial_{\theta_2}$ should be 0, 1 and 1, respectively.

\begin{defn} 
The quantum vector fields $\partial_x, \partial_{\theta_1}, \partial_{\theta_2}:{ O}({\mathbb A}_h^{1|2})\to { O}({\mathbb A}_h^{1|2})$ defined by $(\ref{3.7})$ are called the {\it partial derivatives} of the calculus $(\Omega,{\sf d})$.
\end{defn}
The maps $\partial_x, \partial_{\theta_1}, \partial_{\theta_2}$ are left ${\mathbb K}$-linear with appropriate grade. The family $\{\partial_x, \partial_{\theta_1}, \partial_{\theta_2}\}$ are dual to $\{{\sf d}x, {\sf d}\theta_1, {\sf d}\theta_2\}$.
The partial derivatives of the calculus $(\Omega,{\sf d})$ satisfy the property
\begin{equation*}
\partial_i(x_j) = \delta_{ij}.
\end{equation*}

The next theorem gives the relations between the elements of ${ O}({\mathbb A}_h^{1|2})$ and the partial derivatives.

\begin{thm} 
The partial derivatives $\partial_x$, $\partial_{\theta_1}$ and $\partial_{\theta_2}$, considered as linear mappings of $O({\mathbb A}_h^{1|2})$, and the coordinate functions $x$, $\theta_1$ and $\theta_2$, acting on $O({\mathbb A}_h^{1|2})$ by left multiplication, satisfy the relations
\begin{align} \label{3.8}
\partial_x x & = 1 + x\partial_x + h\, x \partial_{\theta_2} + h'\,\theta_2 \partial_x + hh'\, (x \partial_x + \theta_2 \partial_{\theta_2}), \quad
   \partial_x \theta_1 = \theta_1 \partial_x, \nonumber\\
\partial_x \theta_2 &= \theta_2\partial_x - h(x\partial_x + \theta_2\partial_{\theta_2}), \quad \partial_{\theta_1}x = x \partial_{\theta_1}, \quad
   \partial_{\theta_1} {\theta_1} = 1- {\theta_1} \partial_{\theta_1}, \\
\partial_{\theta_1} \theta_2 &= -\theta_2\partial_{\theta_1}, \quad    \partial_{\theta_2} x = x \partial_{\theta_2} + h'\, (x \partial_x + \theta_2 \partial_{\theta_2}), \quad \partial_{\theta_2}\theta_1 = - \theta_1\partial_{\theta_2}, \nonumber\\
\partial_{\theta_2} \theta_2 &= 1 - \theta_2\partial_{\theta_2} - hx\partial_{\theta_2} + h'\theta_2\partial_x +hh'(x\partial_x + \theta_2\partial_{\theta_2}). \nonumber
\end{align}
\end{thm}

\noindent{\it Proof}
We know that the exterior differential {\sf d} can be written in terms of the differentials and partial derivatives as follows:
\begin{align} \label{3.9}
{\sf d}f &= ({\sf d}x \partial_x + {\sf d}\theta_1 \partial_{\theta_1} + {\sf d}\theta_2 \partial_{\theta_2})f
\end{align}
where $f$ is a differentiable function. Then, if we replace $f$ with $xf$ in the left hand side of the equality in (\ref{3.9}), we get
\begin{align*}
{\sf d} (xf) & = {\sf d}x \, f + x \,({\sf d}x \, \partial_x + {\sf d}\theta_1 \, \partial_{\theta_1} + {\sf d}\theta_2 \, \partial_{\theta_1})f \\
& = \{{\sf d}x \, [1 + x\partial_x + h\, x \partial_{\theta_2} + h'\,\theta_2 \partial_x + hh'\, (x \partial_x + \theta_2 \partial_{\theta_2})] +
 {\sf d}\theta_1 \, (x \, \partial_{\theta_1}) \\
& \quad + {\sf d}\theta_2 \, [x \partial_{\theta_2} + h'\, (x \partial_x + \theta_2 \partial_{\theta_2})]\}f.
\end{align*}
On the other hand, the right hand side of the equality in (\ref{3.9}) has the form
\begin{align*}
{\sf d}(xf) &= [{\sf d}x (\partial_x  x) + {\sf d}\theta_1 (\partial_{\theta_1} x) + {\sf d}\theta_2 (\partial_{\theta_2} x)]f.
\end{align*}
Now, by comparing the right hand sides of these two equalities we obtain some relations in (\ref{3.8}). Other relations can be found similarly. \hfill $\square$

The following theorem can be proved by using the fact that ${\sf d}^2 = 0$ and the relations (\ref{3.1}) together with (\ref{3.9}).

\begin{thm} 
The commutation relations among  the partial derivatives are as follows
\begin{align} \label{3.12}
\partial_x \partial_{\theta_1} &= \partial_{\theta_1} \partial_x, \quad \partial_x \partial_{\theta_2} = \partial_{\theta_2} \partial_x + h'\, \partial_x \partial_x, \quad \partial_{\theta_1} \partial_{\theta_2} = - \partial_{\theta_2} \partial_{\theta_1}, \nonumber\\
\partial_{\theta_1}\partial_{\theta_1} &= 0, \quad \partial_{\theta_2}\partial_{\theta_2} = - h'\, \partial_{\theta_2}\partial_x.
\end{align}
\end{thm}

By Theorem 3.4, the first order differential calculus $\Omega({\mathbb A}_h^{1|2})$ is a $\star$-calculus for the $\star$-algebra
$O({\mathbb A}_h^{1|2})$. The involution of $\Omega({\mathbb A}_h^{1|2})$ induces the involution for the partial differential operators.

Let $D$ be the unital algebra with the generators $\partial_x$, $\partial_{\theta_1}$, $\partial_{\theta_2}$ and defining relations (\ref{3.12}). It can be shown that the set of monomials
\begin{equation*}
\{\partial_x^k \partial_{\theta_1}^l \partial_{\theta_2}^m: \, k\in{\mathbb N}_0, l=0,1, m=0,1,2\}
\end{equation*}
forms a vector space basis of $D$.

\begin{thm} 
The algebra $D$ equipped with the involution determined by
\begin{eqnarray} \label{3.13}
\partial_x^\star = -\partial_x - h \partial_{\theta_2}, \quad \partial_{\theta_1}^\star = \partial_{\theta_1}, \quad \partial_{\theta_2}^\star = \partial_{\theta_2}
\end{eqnarray}
becomes a $\star$-algebra.
\end{thm}

\noindent{\it Proof}
We obtain the equalities (\ref{3.13}) from (\ref{3.5}). It remains to be see that the relations (\ref{3.12}) are preserved under the star operation, which can be checked by direct calculations.
\hfill $\square$

\noindent{\bf Note 4.} It can be also seen that the relations (\ref{3.8}) are preserved the star operation.

\section{The contraction procedure for $\Omega({\mathbb A}_h^{1|2})$} \label{4}

We consider the $q$-deformed algebra of functions on the quantum superspace ${\mathbb A}_q^{1|2}$ generated by $X$, $\Theta_1$ and $\Theta_2$ with the relations \cite{manin1}
\begin{equation} \label{4.1}
X\Theta_i = q\Theta_i X, \quad \Theta_i\Theta_j = -q^{i-j}\Theta_j\Theta_i, \qquad (i,j=1,2)
\end{equation}
and the $(p,q)$-deformed algebra of functions on the quantum superspace ${\mathbb A}_{p,q}^{2|1}$ generated by $\Phi$, $Y_1$, $Y_2$ with the relations \cite{sucelik}
 \begin{equation} \label{4.2}
\Phi^2 = 0, \quad \Phi Y_1 = qp^{-1} Y_1 \Phi, \quad \Phi Y_2 = pq Y_2 \Phi, \quad Y_1Y_2 = pq^{-1} \,Y_2Y_1
\end{equation}
where $p,q\in{\mathbb K}-\{0\}$. The quantum de Rham complex $\Omega({\mathbb A}_q^{1|2})$ is generated by the elements of the set $\{X,\Theta_1,\Theta_2,{\sf d}X,{\sf d}\Theta_1,{\sf d}\Theta_2\}$ by adding the nine cross-commutation relations satisfied between the elements of $O({\mathbb A}_q^{1|2})$ and $O({\mathbb A}_{p,q}^{2|1})$, given in the following lemma, to the relations (\ref{4.1}) and (\ref{4.2}).

\begin{lem} \cite{sucelik} 
The ${\mathbb Z}_2$-graded $(p,q)$-commutation relations are satisfied by the generators of $O({\mathbb A}_q^{1|2})$ and their differentials are stated as below
\begin{align} \label{4.3}
X \,\Phi &=  p \,\Phi \,X,  \quad  X \,Y_1 = q \,Y_1 \,X + (p-1) \,\Phi Y_1,  \quad
  X \,Y_2 = pq \,Y_2 \,X, \nonumber \\
\Theta_1 \,\Phi &= - pq^{-1} \,\Phi \Theta_1,  \quad \Theta_1 \,Y_1 = Y_1 \,\Theta_1,  \quad
  \Theta_1 \,Y_2 = pq^{-1} \,Y_2 \Theta_1, \nonumber\\
\Theta_2 \,\Phi &= - q^{-1} \,\Phi \,\Theta_2 + (1-p) \,Y_2 \,X, \quad \Theta_2 \,Y_1 = q \,Y_1 \,\Theta_2 + (1-p) \,Y_2 \,\Theta_1, \nonumber \\
\Theta_2 \,Y_2 &= Y_2 \,\Theta_2.
\end{align}
\end{lem}

\begin{thm} 
The generators of $O({\mathbb A}_h^{1|2})$ and their differentials satisfy the following commutation relations
\begin{align} \label{4.4}
x \varphi &=  (p+q^{-1} hh') \varphi x + h'(y_2 x -q^{-1}\varphi \theta_2),  \quad  x y_1 = q y_1 x + (p-1) \varphi \theta_1,  \nonumber\\
x y_2 &= pq y_2 x + h (q^{-1}h' \varphi \theta_2 - p\varphi x) + h' y_2 \theta_2, \,\, \theta_1 y_1 = y_1 \theta_1, \,\, \theta_1 y_2 = pq^{-1} y_2 \theta_1, \nonumber \\
\theta_2 \,\varphi &= - q^{-1} \,\varphi \,\theta_2 - q^{-1}h \varphi x - q^{-1} h' y_2 \theta_2 + (1-p-q^{-1}hh') \,y_2 \,x, \\
\theta_2 y_1 &= q y_1 \theta_2 + (1-p) y_2 \theta_1, \quad
\theta_2 y_2 = (1- q^{-1}hh') y_2 \theta_2 - h (q^{-1} \varphi \theta_2 + p y_2 x). \nonumber
\end{align}
\end{thm}

\noindent{\it Proof}
The differential {\sf d} is uniquely determined by the conditions
\begin{equation} \label{4.5}
\Phi = {\sf d}X, \quad Y_1 = {\sf d}\Theta_1, \quad Y_2 = {\sf d}\Theta_2, \quad {\sf d}^2 = 0,
\end{equation}
and the ${\mathbb Z}_2$-graded Leibniz rule.

We now introduce (new) coordinates $x$, $\theta_1$ and $\theta_2$ with the change of basis in the coordinates of the $q$-superspace as follows \cite{celik4}:
\begin{equation} \label{4.6}
X = x + \frac{h'}{pq-1} \, \theta_2, \quad \Theta_1 = \theta_1, \quad \Theta_2 = \theta_2 + \frac{h}{q-1} \, x.
\end{equation}
The parameters $h$ and $h'$ are two deformation parameters that will be replaced with $q$ and $p$ in the limits $q\to1$ and $p\to1$. Considering the differential of the generators of $O({\mathbb A}_q^{1|2})$ and differentiating (\ref{4.6}) we have
\begin{equation} \label{4.7}
{\sf d}x = (1-\tilde{h}\tilde{h'}) \, {\sf d}X -  \tilde{h'} \, {\sf d}\Theta_2, \quad {\sf d}\theta_1 = {\sf d}\Theta_1, \quad
  {\sf d}\theta_2 = (1+\tilde{h}\tilde{h'}) \, {\sf d}\Theta_2 - \tilde{h} \, X
\end{equation}
where $\tilde{h}=h/(q-1)$ and $\tilde{h'}=h'/(pq-1)$. Substituting (\ref{4.6}) and (\ref{4.7}) into (\ref{4.3}), we see that the transformed objects satisfy the relations (\ref{4.4}).
\hfill $\square$

\noindent{\bf Note 5.} While the $q$-deformed calculus is constructed in \cite{wz}, the generators of the exterior algebra have been considered as the differentials of the coordinate functions. In \cite{celik4}, the relations satisfied between the generators of the exterior algebra are normally found with the aid of a singular matrix $g$ containing two Grassmann parameters, using the relations (\ref{4.2}). Here, the generators of the exterior algebra have been obtained after a change of basis in the coordinates of the $q$-superspace (see, above (\ref{4.6}) and (\ref{4.7})). Therefore, while de Rham complex $\Omega(\mathbb A_h^{1|2})$ of the quantum superspace ${\mathbb A}_h^{1|2}$ is established, here it has been given a different definition from that one given in \cite{celik4}.

\noindent{\bf Note 6.} Taking the limits $p\to1$ and $q\to1$ in (\ref{4.4}), one can see that the transformed objects satisfy the commutation relations (\ref{3.2}) given in Theorem 3.2.

\noindent{\bf Note 7.} Substituting (\ref{4.7}) into (\ref{4.2}) and taking the limits $p\to1$ and $q\to1$, one obtains the relations (\ref{3.1}).

If we are introduced commutation relations between the coordinates of the quantum superspace and their partial derivatives after obtaining the derivatives $\partial_x$, $\partial_{\theta_1}$ and $\partial_{\theta_2}$ in terms of $\partial_X$, $\partial_{\Theta_1}$ and $\partial_{\Theta_2}$, we can find the relations (\ref{3.8}) given in Theorem 3.6.

\begin{lem} \cite{sucelik} 
The relations between the generators of $O({\mathbb A}_q^{1|2})$ and partial derivatives are as follows:
\begin{align} \label{4.8}
\partial_X X & = 1 + p X \partial_X + (p-1) \Theta_1 \partial_{\Theta_1}, \quad \partial_X \Theta_1 = pq^{-1} \Theta_1 \partial_X, \nonumber \\
\partial_X \Theta_2 &= q^{-1} \Theta_2 \partial_X, \quad \partial_{\Theta_1} X = q X \partial_{\Theta_1}, \quad \partial_{\Theta_1} \Theta_1 = 1 - \Theta_1\partial_{\Theta_1}, \\
\partial_{\Theta_1} \Theta_2 &= -q \Theta_2 \partial_{\Theta_1},  \quad \partial_{\Theta_2} X = pq X \partial_{\Theta_2}, \quad
 \partial_{\Theta_2} \Theta_1 = - pq^{-1} \Theta_1 \partial_{\Theta_2}, \nonumber \\
\partial_{\Theta_2} \Theta_2 &= 1 - \Theta_2 \partial_{\Theta_2} + (p-1) (X \partial_X + \Theta_1 \partial_{\Theta_1}).\nonumber
\end{align}
\end{lem}

\begin{thm} 
The relations of the generators of $O({\mathbb A}_h^{1|2})$ with partial derivatives are as follows:
\begin{align} \label{4.9}
\partial_x \, x & = 1 + (p+q^{-1}hh') x \, \partial_x + ph \, x \, \partial_{\theta_2} + q^{-1}h' \theta_2 \partial_x + q^{-1}hh' \theta_2 \partial_{\theta_2} \nonumber \\ & \quad + (p-1) \theta_1 \partial_{\theta_1}, \nonumber \\
\partial_x \, \theta_1 &= pq^{-1} \theta_1 \, \partial_x, \quad \partial_x \, \theta_2 = q^{-1} \theta_2 \, \partial_x - q^{-1} h (x \partial_x +  \theta_2 \partial_{\theta_2}), \\
\partial_{\theta_1} x & = q x \partial_{\theta_1}, \quad \partial_{\theta_1} \theta_1 = 1- \theta_1 \partial_{\theta_1}, \quad \partial_{\theta_1} \theta_2 = -q \theta_2 \partial_{\theta_1}, \nonumber \\
\partial_{\theta_2} x & = pq x \partial_{\theta_2} + h' (x \partial_x + \theta_2 \partial_{\theta_2}), \quad
 \partial_{\theta_2} \theta_1 = -pq^{-1} \theta_1 \partial_{\theta_2}, \nonumber \\
\partial_{\theta_2} \theta_2 &= 1 - (1-q^{-1}hh') \theta_2 \partial_{\theta_2} + ph x \partial_{\theta_2} + q^{-1}h' \theta_2 \partial_x + (p-1+q^{-1}hh') x \partial_x \nonumber \\
& \quad + (p-1) \theta_1 \partial_{\theta_1}. \nonumber
\end{align}
\end{thm}

\noindent{\it Proof}
From the equalities in (\ref{4.6}), we have
\begin{align} \label{4.10}
\partial_X &= \left(1-\tfrac{hh'}{(q-1)(pq-1)}\right) \partial_x - \tfrac{h}{q-1} \, \partial_{\theta_2}, \quad \partial_{\Theta_1} = \partial_{\theta_1}, \nonumber \\
\partial_{\Theta_2} &= \left(1+\tfrac{hh'}{(q-1)(pq-1)}\right) \partial_{\theta_2} + \tfrac{h'}{pq-1} \partial_x.
\end{align}
Substituting (\ref{4.6}) and (\ref{4.10}) into (\ref{4.8}), we obtain the relations (\ref{4.9}).
\hfill $\square$

\noindent{\bf Note 8.} We know that, in the $q$-deformed case, the exterior differential can be written in terms of the differentials and partial derivatives as follows:
\begin{equation*}
{\sf d}'f = ({\sf d}X \, \partial_X + {\sf d}\Theta_1 \, \partial_{\Theta_1} + {\sf d}\Theta_2 \, \partial_{\Theta_2})f
\end{equation*}
where $f$ is a differentiable function. If (\ref{4.7}) and (\ref{4.10}) are considered, then the differential {\sf d} preserves its form
\begin{equation*}
{\sf d}f = ({\sf d}x \, \partial_x + {\sf d}\theta_1 \, \partial_{\theta_1} + {\sf d}\theta_2 \, \partial_{\theta_2})f.
\end{equation*}


\end{document}